\documentclass{amsart}

\usepackage{amsthm,amsmath,amssymb,eepic,graphicx}
\usepackage[usenames]{color}

\title
[CMC $1$ surfaces in de Sitter $3$-space]
{Spacelike mean curvature $1$ surfaces of genus $1$ \\ 
 with two ends in de Sitter $3$-Space}
\author{Shoichi FUJIMORI}
\date{May 19, 2006}
\address{Faculty of Mathematics, Kyushu University, Fukuoka 812-8581, Japan}
\email{fujimori@math.kyushu-u.ac.jp}
\thanks{2000 {\it Mathematics Subject Classification.} 53A10, 53B30. \\
\indent
{\it Key words and phrases.} de Sitter $3$-space, 
                             spacelike CMC $1$ surface, 
                             genus $1$ surface.}
\def\transpose#1{\mathord{\mathopen{{\vphantom{#1}}^t}#1}}

\DeclareFontFamily{U}{rsfs}{\skewchar\font"7F}
\DeclareFontShape{U}{rsfs}{m}{n}{
	<-6> rsfs5
	<6-8> rsfs7
	<8-> rsfs10
	}{}
\DeclareMathAlphabet{\mathscr}{U}{rsfs}{m}{n}
 \newtheorem{theorem}{Theorem}[section]
 \newtheorem*{theorem*}{Theorem}
 \newtheorem{proposition}[theorem]{Proposition}
 \newtheorem{corollary}[theorem]{Corollary}
 \newtheorem{lemma}[theorem]{Lemma}

\theoremstyle{definition}
 \newtheorem{definition}[theorem]{Definition}
 \newtheorem{example}[theorem]{Example}
\theoremstyle{remark}
 \newtheorem{remark}[theorem]{Remark}
 \newtheorem*{remark*}{Remark}

\numberwithin{equation}{section}
\begin{document} %%%%%%%%%%%%%%%%%%%%%%%%%%%%%%%%%%%%%%%%%%%%%%%%%%%%%%%

\begin{abstract}
We give a mathematical foundation for, and numerical 
demonstration of, the existence of mean curvature $1$ surfaces of genus $1$ 
with either two elliptic ends or two hyperbolic ends in de Sitter $3$-space. 
An end of a mean 
curvature $1$ surface is an ``elliptic end'' (resp. a ``hyperbolic end'') if 
the monodromy matrix at the end is diagonalizable with eigenvalues in the unit 
circle (resp. in the reals).  
Although the existence of the surfaces is numerical, 
the types of ends are mathematically determined.  
\end{abstract}

\maketitle

%\vspace{0.5\baselineskip}
%
\section*{Introduction} %%%%%%%%%%%%%%%%%%%%%%%%%%%%%%%%%%%%%%%%%%%%%%%%
%%%%%%%%%%%%%%%%%%%%%%%%%%%%%%%%%%%%%%%%%%%%%%%%%%%%%%%%%%%%%%%%%%%%%%%%

The global theories of minimal surfaces in Euclidean $3$-space $\mathbb{R}^3$ 
and constant mean curvature (CMC) $1$ surfaces in hyperbolic $3$-space 
$\mathbb{H}^3$ are well understood, as they possess representation formulas 
using meromorphic functions and so benefit from the theory of complex 
analysis. 

In contrast to this, the global theory of spacelike maximal surfaces in 
Minkowski $3$-space $\mathbb{R}^3_1$ and spacelike CMC $1$ surfaces 
in de Sitter $3$-space $\mathbb{S}^3_1$ are not well explored yet, 
even though they possess similar representation formulas. 
This is perhaps because the only complete spacelike maximal immersions in 
$\mathbb{R}^3_1$ and spacelike CMC $1$ immersions in $\mathbb{S}^3_1$ 
are flat and totally umbilic. 
So to have an interesting global theory about these surfaces, we need to 
consider a wider class of surfaces than just complete and immersed ones. 

Recently, Umehara and Yamada defined such a category of spacelike maximal 
surfaces with certain kinds of singularities and named them ``maxfaces'' 
\cite{UYmax}.  Then they constructed numerous examples by a transferring 
method from minimal surfaces in $\mathbb{R}^3$. 
Furthermore, Kim and Yang discovered an interesting example of a maxface, 
which has genus $1$ with two embedded ends, even though there does not exist 
such an example as a complete minimal immersion in $\mathbb{R}^3$ \cite{KY}. 
In addition, Fern\`andez and L\'opez and Souam have investigated maximal 
surfaces with conical singularities \cite{FLS1, FLS2}. 

The author defined spacelike CMC $1$ surfaces with certain kinds of 
singularities as an analogue of maxfaces, naming them ``CMC $1$ faces'', 
and constructed many examples by transferring from reducible CMC $1$ surfaces 
in $\mathbb{H}^3$ \cite{F}. 
Also, Lee and Yang investigated spacelike CMC $1$ surfaces of genus zero 
with two and three ends \cite{LY}.  
However, every surface constructed in \cite{F} and \cite{LY} was 
topologically a sphere with finitely many points removed. 
Given all of this, it is natural to consider 
whether or not there exist examples with positive genus. 

For CMC $1$ immersions in $\mathbb{H}^3$, Rossman and Sato constructed genus 
$1$ catenoid cousins by a numerical method \cite{RS}. 
Here we will similarly construct genus $1$ ``catenoids'' using a modification 
of their method;  
that is, we show the following numerical result (Example \ref{ex:g1cat}):
\begin{quote}
There exist one-parameter families of weakly-complete CMC $1$ faces 
of genus $1$ with two elliptic or two hyperbolic ends which satisfy 
equality in the Osserman-type inequality. 
\end{quote}

The Osserman inequality for complete minimal immersions in $\mathbb{R}^3$ says 
that twice the degree of the Gauss map is greater than or equal to the number 
of ends minus the Euler characteristic of the surface, with equality holding 
if and only if all the ends are embedded.  
An analogous Osserman-type inequality for CMC $1$ faces in $\mathbb{S}^3_1$ 
was shown in \cite{F}, in the case that the ends are complete and elliptic.  
The examples here satisfy  equality in the Osserman-type inequality, 
even though some of them do not have elliptic ends. 
(We define elliptic and hyperbolic and parabolic ends in 
 Section \ref{sec:prelim}.)
Osserman-type inequalities for CMC $1$ immersions in $\mathbb{H}^3$ and 
maxfaces in $\mathbb{R}^3_1$ can be found in \cite{UY1, UY5} and \cite{UYmax}. 

For weakly-complete CMC $1$ faces, 
the behavior of ends is investigated in \cite{FRUYY}. 
In addition, criteria for the singularities are given in \cite{FSUY}.

The author would like to thank Professors Wayne Rossman, Masaaki Umehara, 
Kotaro Yamada and Seong-Deog Yang for their valuable comments and suggestions. 

\section{Preliminaries} %%%%%%%%%%%%%%%%%%%%%%%%%%%%%%%%%%%%%%%%%%%%%%%%
\label{sec:prelim} %%%%%%%%%%%%%%%%%%%%%%%%%%%%%%%%%%%%%%%%%%%%%%%%%%%%%

\subsection{de Sitter $3$-space} %%%%%%%%%%%%%%%%%%%%%%%%%%%%%%%%%%%%%

Let $\mathbb{R}^4_1$ be the $4$-dimensional Lorentz space with the Lorentz 
metric 
\[ 
\langle (x_0,x_1,x_2,x_3),(y_0,y_1,y_2,y_3)\rangle =
-x_0y_0+x_1y_1+x_2y_2+x_3y_3. 
\]
Then de Sitter $3$-space is 
\[
\mathbb{S}^3_1=\mathbb{S}^3_1(1)=
\{(x_0,x_1,x_2,x_3)\in\mathbb{R}^4_1 \, | \, -x_0^2+x_1^2+x_2^2+x_3^2=1\},
\]
with metric induced from $\mathbb{R}^4_1$. 
$\mathbb{S}^3_1$ is a simply-connected %complete 
$3$-dimensional Lorentzian manifold with constant sectional curvature $1$. 
We can consider $\mathbb{R}^4_1$ to be the $2 \times 2$ self-adjoint matrices 
($X^* = X$, where $X^*=\transpose{\overline{X}}$, and $\transpose{X}$ denotes 
the transpose of $X$) by the identification 
\[
\mathbb{R}^4_1\ni X=(x_0,x_1,x_2,x_3)\leftrightarrow 
X=\sum_{k=0}^3 x_k e_k 
 =\begin{pmatrix} x_0+x_3           & x_1+i x_2 \\ 
                  x_1-i x_2 & x_0-x_3 \end{pmatrix},
\]
where
\[
e_0=\begin{pmatrix}1&0\\0&1\end{pmatrix},\quad
e_1=\begin{pmatrix}0&1\\1&0\end{pmatrix},\quad
e_2=\begin{pmatrix}0&i\\-i&0\end{pmatrix},\quad
e_3=\begin{pmatrix}1&0\\0&-1\end{pmatrix}.
\]
Then $\mathbb{S}^3_1$ is 
\[ 
\mathbb{S}^3_1=\{X\,|\,X^*=X\,,\det X=-1\}
              =\{Fe_3F^*\,|\,F\in SL(2,\mathbb{C})\} 
\] 
with the metric
\[
 \langle X,Y\rangle 
=-\frac{1}{2}\mathrm{trace}\left(Xe_2(\transpose{Y})e_2\right) . 
\]
In particular, $\langle X,X\rangle =-\det X$. 
An immersion in $\mathbb{S}^3_1$ is called {\em spacelike} if the induced 
metric on the immersed surface is positive definite.  

\subsection{CMC $1$ faces} %%%%%%%%%%%%%%%%%%%%%%%%%%%%%%%%%%%%%%%%%%%

Aiyama and Akutagawa gave a local Weierstrass-type representation formula for 
spacelike immersions of constant mean curvature (CMC) $1$ in 
$\mathbb{S}^3_1$ \cite{AA}. 
However, for complete spacelike CMC $1$ immersions in $\mathbb{S}^3_1$, 
the only ones that exist are totally umbilic \cite{Ak, R}. 
So we must enlarge the class of surfaces we consider to include 
non-immersions, in order to have an interesting theory: 

\begin{definition}
Let $M$ be an oriented $2$-manifold. 
A $C^\infty$-map $f:M\to\mathbb{S}^3_1$ is called a {\em CMC $1$ face} 
\cite{F} if 
\begin{enumerate}
\item there exists an open dense subset $W\subset M$ such that 
$f|_W$ is a spacelike CMC $1$ immersion, 
\item for any singular point $p$ 
(that is, a point where the induced metric degenerates), 
there exists a $C^1$-differentiable function 
$\lambda :U\cap W\to\mathbb{R}^+$, where $U$ is a neighborhood of $p$, 
such that $\lambda ds^2$ extends to a $C^1$-differentiable 
Riemannian metric on $U$, and 
\item $df(p)\ne 0$ for any $p\in M$. 
\end{enumerate}
\end{definition}

It is known that the $2$-manifold $M$ on which a CMC $1$ face 
$f:M\to\mathbb{S}^3_1$ is defined always has a complex structure \cite{F}.  
So we will treat $M$ as a Riemann surface. 

The representation formula of Aiyama-Akutagawa can be extended 
to CMC $1$ faces as follows:

\begin{theorem}\label{th:AA-rep} \cite{F}
Let $M$ be a Riemann surface with a base point $z_0\in M$. Let $G$ be a 
meromorphic function and $Q$ a holomorphic $2$-differential on $M$ such that 
\begin{equation}\label{eq:dshat^2}
ds_\#^2=(1+|G|^2)^2\frac{Q}{dG}\overline{\left(\frac{Q}{dG}\right)}
\end{equation} 
is a Riemannian metric on $M$. 
Choose the holomorphic immersion $F=(F_{jk})$ defined on the universal cover 
$\widetilde M$ of $M$ into $SL(2,\mathbb{C})$ 
so that $F(z_0)=e_0$ and $F$ satisfies 
\begin{equation}\label{eq:F^-1dF}
dF\cdot F^{-1}=\alpha\quad\text{where}\quad
  \alpha=\begin{pmatrix}G&-G^2\\1&-G\end{pmatrix}\frac{Q}{dG}.
\end{equation}
Then $f:\widetilde M\to\mathbb{S}^3_1$ defined by 
\begin{equation}\label{eq:f=Fe_3F^*}
f=Fe_3F^*
\end{equation}
is a CMC $1$ face that is conformal away from its singularities. The induced 
metric $ds^2$ on $M$ and the second fundamental form $h$ are given as follows:
\begin{equation}\label{eq:ds^2-h-G}
ds^2=(1-|g|^2)^2\frac{Q}{dg}\overline{\left(\frac{Q}{dg}\right)} ,\quad 
   h=Q+\overline{Q}+ds^2,
\end{equation}
where $g$ is defined as the multi-valued function 
$-dF_{12}/dF_{11}=-dF_{22}/dF_{21}$ on $M$. 
Moreover, $G$ is the hyperbolic Gauss map of $f$ and 
$Q$ is the Hopf differential of $f$. 
The singularities of the CMC $1$ face occur at points where $|g|=1$. 

Conversely, 
let $M$ be a Riemann surface and $f:M\to\mathbb{S}^3_1$ a CMC $1$ face. 
Then there exists a meromorphic function $G$ 
and holomorphic $2$-differential $Q$ on $M$ such that $ds_\#^2$ is a 
Riemannian metric on $M$, and such that \eqref{eq:f=Fe_3F^*} holds, 
where $F:\widetilde{M}\to SL(2,\mathbb{C})$ is an immersion which satisfies 
\eqref{eq:F^-1dF}. 
\end{theorem}

\begin{remark}\label{rm:AArep}
We make the following remarks about Theorem \ref{th:AA-rep}:
\begin{enumerate}
\item Following the terminology of Umehara and Yamada, $g$ is called the 
{\em secondary} Gauss map. We call $(G,Q)$ the {\em Weierstrass data}, 
and $F$ the {\em holomorphic null lift} of $f$. 
Also, $ds_\#^2$ defined as in \eqref{eq:dshat^2} is called the 
{\em lift metric} of $f$. 
\item For a regular point, the unit normal vector 
$N$ of $f$ is given by 
\begin{equation}\label{eq:unitnormal}
       N=\frac{1}{|g|^2-1}(F\nu )(F\nu )^*, \quad\text{where}\quad
    \nu =\begin{pmatrix}1&g\\\bar g&1\end{pmatrix}, 
\end{equation}
which is a future pointing (resp. past pointing) vector 
if and only if $|g|>1$ (resp. $|g|<1$). 
We also remark that $N$ is a unit {\em timelike} vector, that is, 
$\langle N,N\rangle =-1$. 
\item When $|g|>1$ (resp. $|g|<1$), 
the hyperbolic Gauss map has the following geometric meaning: 
Let $\mathbb{S}^2_\infty\cong\mathbb{C}\cup\{\infty\}$ be the future 
(resp. past) pointing ideal boundary of $\mathbb{S}^3_1$. 
Let $\gamma_z$ be the geodesic ray starting at $f(z)$ in $\mathbb{S}^3_1$ 
with the velocity vector $N(z)$ at $f(z)$. 
Then $G(z)$ is the point in $\mathbb{S}^2_\infty$ determined by the asymptotic 
class of $\gamma_z$. 
See \cite{B, UY1, FRUYY}. 
\item\label{rm:2Q=Sg-SG} 
By Equation (2.6) in \cite{UY1}, $G$ and $g$ and $Q$ have the following 
relation: 
\begin{equation}\label{eq:2Q=Sg-SG}
2Q=S(g)-S(G),
\end{equation}
where $S(g)=S_z(g)dz^2$ and
\[
S_z(g)=\left(\frac{g''}{g'}\right)'-\frac{1}{2}\left(\frac{g''}{g'}\right)^2
\qquad \left({}'=\frac{d}{dz}\right)
\]
is the Schwarzian derivative of $g$. 
\item 
For a CMC $1$ face $f$, if we find both the hyperbolic Gauss map $G$ and the 
secondary Gauss map $g$, we can explicitly find the holomorphic null lift $F$, 
by using the so-called Small formula:
\[
F=\begin{pmatrix}G\dfrac{da}{dG}-a&G\dfrac{db}{dG}-b\\[6pt]
                 \phantom{G}\dfrac{da}{dG}\phantom{-a}&
                 \phantom{G}\dfrac{db}{dG}\phantom{-b}\end{pmatrix}, \quad
a=\sqrt{\dfrac{dG}{dg}},\quad b=-ga. 
\]
See \cite{S, KUYsmall}. 
\end{enumerate}
\end{remark}

\subsection{Closing conditions for CMC $1$ faces} %%%%%%%%%%%%%%%%%%%%

\begin{definition}\label{df:comp-fin}
Let $M$ be a Riemann surface and $f:M\to\mathbb{S}^3_1$ a CMC $1$ face. 
Set $ds^2=f^*(ds^2_{\mathbb{S}^3_1})$. 
\begin{enumerate}
\item 
$f$ is {\em complete} (resp. of {\em finite type}) if there exists a compact 
set $C$ and a symmetric $(0,2)$-tensor $T$ on $M$ such that $T$ vanishes on 
$M\setminus C$ and $ds^2+T$ is a complete (resp. finite total curvature) 
Riemannian metric. 
\item $f$ is {\em weakly-complete} 
(resp. of {\em weakly finite total curvature}) if the 
lift metric $ds_\#^2$ defined as in \eqref{eq:dshat^2} is a complete 
(resp. finite total curvature) Riemannian metric (\cite{FRUYY}). 
\end{enumerate}
\end{definition}

\begin{remark}
If a CMC $1$ face $f$ is complete and of finite type, then $f$ is 
weakly-complete and of weakly finite total curvature \cite{F}. 
But the converse is not true.  See \cite{FRUYY}. 
\end{remark}

Let $f:M\to\mathbb{S}^3_1$ be a CMC $1$ face defined on a Riemann surface $M$ 
biholomorphic to a compact Riemann surface with finitely many points removed. 
The ends of $f$ correspond to the removed points. 
Let $\varrho :\widetilde{M}\to M$ be the universal cover of $M$, and 
$F:\widetilde{M}\to SL(2,\mathbb{C})$ a holomorphic null lift of $f$.  
We fix a point $z_0\in M$.  
Let $\gamma :[0,1]\to M$ be a loop so that $\gamma (0)=\gamma (1)=z_0$. 
Then there exists a unique deck transformation $\tau$ of $\widetilde{M}$ 
associated to the homotopy class of $\gamma$.  We define the monodromy matrix 
$\Phi_\gamma$ of $F$ with respect to $\gamma$ by 
\[
F\circ\tau =F\Phi_\gamma.
\]
If a loop $\gamma$ lies in a small neighborhood of an end of $f$ and wraps 
once (has winding number $\pm 1$) about the end, then $\Phi_\gamma$ can be 
regarded as the monodromy matrix about the end.  
We give the following definition:

\begin{definition}
$F$ satisfies the {\em $SU(1,1)$ condition} if $\Phi_\gamma\in SU(1,1)$ 
for any loop $\gamma$ in $M$. 
\end{definition}

\begin{remark}
Note that $f$ is well-defined on $M$ if and only if $F$ satisfies the 
$SU(1,1)$ condition.  
\end{remark}

Now we assume that $f$ is well-defined on $M$, 
so $F$ does satisfy the $SU(1,1)$ condition.
Then $\Phi_\gamma$ is conjugate to either 
\begin{equation}\label{eq:e-h-p}
\mathcal{E}=\begin{pmatrix}e^{i\theta}&0\\0&e^{-i\theta}\end{pmatrix}
\quad\text{or}\quad
\mathcal{H}=\pm\begin{pmatrix}\cosh s&\sinh s\\\sinh s&\cosh s\end{pmatrix}
\quad\text{or}\quad
\mathcal{P}=\pm\begin{pmatrix}1\pm i&1\\1&1\mp i\end{pmatrix}
\end{equation}
for $\theta\in [0,2\pi)$, $s\in\mathbb{R}\setminus\{0\}$.  

\begin{definition}\label{df:e-h-p}
Let $f:M\to\mathbb{S}^3_1$ be a CMC $1$ face with 
holomorphic null lift $F$. 
An end of $f$ is called an {\em elliptic end} or {\em hyperbolic end} or 
{\em parabolic end} if the monodromy about the end is 
conjugate to $\mathcal{E}$ or $\mathcal{H}$ or $\mathcal{P}$ 
in $SU(1,1)$, respectively.  
An end of $f$ is called {\em regular} if the hyperbolic Gauss map extends 
meromorphically to the end.  
\end{definition}

\begin{remark}
When $F$ changes to $F\Phi$ by a deck transformation, 
the secondary Gauss map $g$ changes to 
\[
\Phi^{-1}\star g
:=\frac{\phantom{-}\Phi_{22}g-\Phi_{12}}{-\Phi_{21}g+\Phi_{11}}, 
\quad\text{where}\quad
\Phi =\begin{pmatrix}\Phi_{11}&\Phi_{12}\\\Phi_{21}&\Phi_{22}\end{pmatrix}. 
\]
So we can consider the $SU(1,1)$ condition as an $SU(1,1)$ period condition 
for $g$. 
\end{remark}

\subsection{Osserman-type inequality} %%%%%%%%%%%%%%%%%%%%%%%%%%%%%%%%

A complete CMC $1$ face of finite type with elliptic ends has the following 
property: 

\begin{theorem}[Osserman-type inequality \cite{F}]\label{th:Oss-ineq}
Let $f:M\to\mathbb{S}^3_1$ be a complete CMC $1$ face of finite type with $n$ 
elliptic ends and no other ends.  
Then $M$ is biholomorphic to $\overline{M}\setminus\{p_1,\dots ,p_n\}$, 
where $\overline{M}$ is a compact Riemann surface. 
Let $G$ be its hyperbolic Gauss map.  Then the following inequality holds:
\begin{equation}\label{eq:Oss-ineq}
2\deg (G)\geq -\chi (M)+n, 
\end{equation}
where $\deg (G)$ is the mapping degree of $G$ $($if $G$ has essential 
singularities, then we define $\deg (G)=\infty)$. 
Furthermore, equality holds if and only if each end is regular and embedded. 
\end{theorem}

\subsection{The hollow ball model} %%%%%%%%%%%%%%%%%%%%%%%%%%%%%%%%%%%

To visualize CMC $1$ faces, we use the {\em hollow ball model} of 
$\mathbb{S}^3_1$, as in \cite{LY}. 
For any point
\[
\begin{pmatrix}x_0+x_3&x_1+ix_2\\x_1-ix_2&x_0-x_3\end{pmatrix}\leftrightarrow
(x_0,x_1,x_2,x_3)\in\mathbb{S}^3_1,
\]
define
\[
y_k=\frac{e^{\arctan x_0}}{\sqrt{1+x_0^2}}x_k,
\qquad k=1,2,3.
\]
Then $e^{-\pi}<y_1^2+y_2^2+y_3^2<e^\pi$. 
The identification $(x_0,x_1,x_2,x_3)\leftrightarrow (y_1,y_2,y_3)$ is then a 
bijection from $\mathbb{S}^3_1$ to the hollow ball 
\[
\mathscr{H}=\{(y_1,y_2,y_3)\in\mathbb{R}^3\,|\,
              e^{-\pi}<y_1^2+y_2^2+y_3^2<e^\pi\}. 
\]
So $\mathbb{S}^3_1$ is identified with the hollow ball $\mathscr{H}$, and we 
show the graphics in this paper using this identification to $\mathscr{H}$. 

\section{CMC $1$ faces of genus $1$} %%%%%%%%%%%%%%%%%%%%%%%%%%%%%%%%%%%
\label{sec:period} %%%%%%%%%%%%%%%%%%%%%%%%%%%%%%%%%%%%%%%%%%%%%%%%%%%%%

Consider the hyperelliptic Riemann surface 
\begin{equation}\label{eq:riemsurf}
M=\left\{(z,w)\in (\mathbb{C}\cup\{\infty\})^2\,\left|\,
         w^2=\frac{(z+1)(z-a)}{(z-1)(z+a)}\right.\right\}
  \setminus\biggl\{(\infty,1), (\infty,-1)\biggr\},
\end{equation}
where $a>1$. Then $M$ is a twice punctured torus. 
Define 
\begin{equation}\label{eq:W-data}
G=w,\qquad Q=\frac{cdzdw}{w} 
\end{equation}
for $c\in\mathbb{R}\setminus\{0\}$. 
Then $ds^2_\#$ defined as in Equation \eqref{eq:dshat^2} gives a Riemannian 
metric on $M$. 
So $(G,Q)$ are the Weierstrass data for a 
genus $1$ catenoid.  Let $F(z,w)\in SL(2,\mathbb{C})$ be the solution of 
Equation \eqref{eq:F^-1dF} with initial condition $F(0,1)=e_0$. 
Then $f=Fe_3F^*$ is a CMC $1$ face in $\mathbb{S}^3_1$, and this 
CMC $1$ face is defined on the universal cover $\widetilde{M}$ of $M$. 

\begin{figure}[htbp] %%%%%%%%%%%%%%%%%%%%%%%%%%%%%%%%%%%%%%%%%%%%%%%%%%%
\begin{center}
\unitlength=1.0pt
\begin{picture}(300,200)
\put(150,100){\circle{170}}
\put(150,100){\circle{70}}
\bezier100(150,185)(130,160)(150,135)
\bezier40(150,185)(170,160)(150,135)
\bezier100(150,65)(130,40)(150,15)
\bezier40(150,65)(170,40)(150,15)
\put(150,185){\circle*{2}}%point at a
\put(150,135){\circle*{2}}%point at 1
\put(150,65){\circle*{2}}%point at -1
\put(150,15){\circle*{2}}%point at -a
\put(65,100){\circle*{2}}%point at (infty,1)
\put(115,100){\circle*{2}}%point at (0,1)
\put(185,100){\circle*{2}}%point at (0,-1)
\put(235,100){\circle*{2}}%point at (infty,-1)
\put(150,192){\makebox(0,0)[cc]{$z=a$}}
\put(150,128){\makebox(0,0)[cc]{$z=1$}}
\put(150,72){\makebox(0,0)[cc]{$z=-1$}}
\put(150,8){\makebox(0,0)[cc]{$z=-a$}}
\put(63,100){\makebox(0,0)[rc]{$(z,w)=(\infty,1)$}}
\put(117,100){\makebox(0,0)[lc]{$(0,1)$}}
\put(183,100){\makebox(0,0)[rc]{$(0,-1)$}}
\put(237,100){\makebox(0,0)[lc]{$(z,w)=(\infty,-1)$}}
\put(150,100){\arc{90}{-2.5}{2.5}}%gamma_1
\bezier70(115,100)(105,114)(114,127)%gamma_1
\bezier70(115,100)(105,86)(114,73)%gamma_1
\put(195,100){\vector(0,-1){0}}%gamma_1
\put(197,100){\makebox(0,0)[lb]{$\gamma_1$}}
\bezier150(115,100)(80,137)(108,174)%gamma_2
\bezier60(115,100)(120,137)(108,174)%gamma_2
\put(96,137){\vector(0,-1){0}}%gamma_2
\put(94,137){\makebox(0,0)[rb]{$\gamma_2$}}
\bezier150(115,100)(80,63)(108,26)%gamma_3
\bezier60(115,100)(120,63)(108,26)%gamma_3
\put(96,63){\vector(0,-1){0}}%gamma_3
\put(94,63){\makebox(0,0)[rt]{$\gamma_3$}}
\end{picture}
\end{center}
\caption{The Riemann surface $M$.  This picture indicates how the loops 
         $\gamma_1$, $\gamma_2$ and $\gamma_3$ lie in $M$.}
\label{fg:torus}
\end{figure}
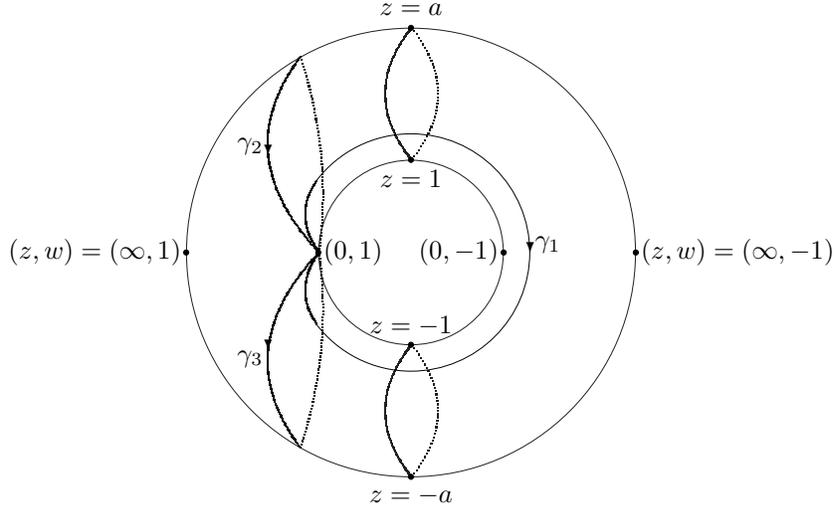

We do not yet know that $f$ is well-defined on $M$ itself. 
For this to happen, $F$ must satisfy the $SU(1,1)$ condition. 
We satisfy the $SU(1,1)$ condition by changing the initial condition $F(0,1)$. 
It is enough to check the $SU(1,1)$ condition on the following three loops, 
since they generate the fundamental group of $M$ 
(see Figures \ref{fg:loops} and \ref{fg:torus}): 

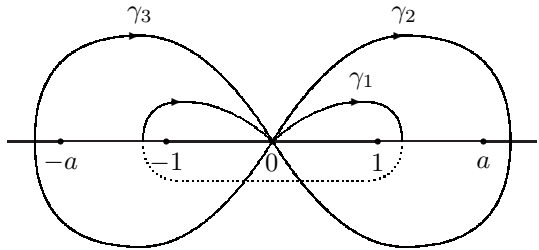
\begin{figure}[htbp] %%%%%%%%%%%%%%%%%%%%%%%%%%%%%%%%%%%%%%%%%%%%%%%%%%%
\begin{center}
\unitlength=1.0pt
\begin{picture}(300,100)
\put(50,50){\line(1,0){200}}%real axis
\thicklines
%\put(150,0){\line(0,1){100}}%imaginary axis
\put(50,50){\line(1,0){20}}%real axis
\put(110,50){\line(1,0){80}}%real axis
\put(230,50){\line(1,0){20}}%real axis
\thinlines
\put(70,50){\circle*{2}}%point at -a
\put(110,50){\circle*{2}}%point at -1
\put(150,50){\circle*{2}}%point at 0
\put(190,50){\circle*{2}}%point at 1
\put(230,50){\circle*{2}}%point at a
\put(70,42){\makebox(0,0)[cc]{$-a$}}
\put(110,42){\makebox(0,0)[cc]{$-1$}}
\put(150,42){\makebox(0,0)[cc]{$0$}}
\put(190,42){\makebox(0,0)[cc]{$1$}}
\put(230,42){\makebox(0,0)[cc]{$a$}}
\put(184,73){\makebox(0,0)[cc]{$\gamma_1$}}
\put(184,65){\vector(1,0){0}}
\put(116,65){\vector(1,0){0}}
\bezier57(150,50)(167,65)(184,65)%gamma1
\bezier40(184,65)(199,65)(199,50)%gamma1
\bezier12(199,50)(199,35)(184,35)%gamma1
\bezier36(184,35)(150,35)(116,35)%gamma1
\bezier12(116,35)(101,35)(101,50)%gamma1
\bezier40(101,50)(101,65)(116,65)%gamma1
\bezier57(116,65)(133,65)(150,50)%gamma1
\put(200,98){\makebox(0,0)[cc]{$\gamma_2$}}
\put(200,90){\vector(1,0){0}}
\bezier150(150,50)(175,90)(200,90)%gamma2
\bezier150(200,90)(240,90)(240,50)%gamma2
\bezier150(240,50)(240,10)(200,10)%gamma2
\bezier150(200,10)(175,10)(150,50)%gamma2
\put(100,98){\makebox(0,0)[cc]{$\gamma_3$}}
\put(100,90){\vector(1,0){0}}
\bezier150(150,50)(125,10)(100,10)%gamma3
\bezier150(100,10)(60,10)(60,50)%gamma3
\bezier150(60,50)(60,90)(100,90)%gamma3
\bezier150(100,90)(125,90)(150,50)%gamma3
\end{picture}
\end{center}
\caption{Projection to the $z$-plane of the loops $\gamma_1$, $\gamma_2$, 
         and $\gamma_3$ which generate the fundamental group of $M$.}
\label{fg:loops}
\end{figure} %%%%%%%%%%%%%%%%%%%%%%%%%%%%%%%%%%%%%%%%%%%%%%%%%%%%%%%%%%%

 %%%%%%%%%%%%%%%%%%%%%%%%%%%%%%%%%%%%%%%%%%%%%%%%%%%%%%%%%%%

\begin{itemize}
\item 
The curve $\gamma_1:[0,1]\to M$ starts at $\gamma_1(0)=(0,1)\in M$. 
Its first portion has $z$ coordinate in the first quadrant of the $z$ 
plane and ends at a point $(z,w)$ where $z\in\mathbb{R}$ and $1<z<a$.  
Its second portion starts at $(z,w)$ and ends at $(0,-1)$ and has $z$ 
coordinate in the fourth quadrant.  
Its third portion starts at $(0,-1)$ and ends at $(-z,1/w)$ and 
has $z$ coordinate in the third quadrant.  
Its fourth and last portion starts at $(-z,1/w)$ and returns to the base point 
$\gamma_1(1)=(0,1)$ and has $z$ coordinate in the second quadrant.  
\item
The curve $\gamma_2:[0,1]\to M$ starts at $\gamma_2(0)=(0,1)$.  
Its first portion has $z$ coordinate in the first quadrant and ends at a 
point $(z,w)$ where $z\in\mathbb{R}$ and $z>a$.  
Its second and last portion starts at $(z,w)$ and returns to 
$\gamma_2(1)=(0,1)$ and has $z$ coordinate in the fourth quadrant.  
\item
The curve $\gamma_3:[0,1]\to M$ starts at $\gamma_3(0)=(0,1)$.  
Its first portion has $z$ coordinate in the third quadrant and ends at a 
point $(z,w)$ where $z\in\mathbb{R}$ and $z<-a$.  
Its second and last portion starts at $(z,w)$ and returns to 
$\gamma_3(1)=(0,1)$ and has $z$ coordinate in the second quadrant.  
\end{itemize}

Consider the symmetries
\[
\begin{array}{ll}
  \phi_1(z,w)=(\bar z,\bar w),\quad
& \phi_2(z,w)=(-z,1/w), \\
  \phi_3(z,w)=(-\bar z,1/\bar w),\quad
& \phi_4(z,w)=(\bar z,-\bar w)
\end{array}
\]
on $M$.  Then by the same argument as Lemmas 5.1 and 5.2 in \cite{RS}, 
we have the following:

\begin{lemma}[{\cite[Lemmas 5.1 and 5.2]{RS}}]
If 
\[
F(z,w)=\begin{pmatrix}A&B\\C&D\end{pmatrix},
\]
then
\begin{equation}\label{eq:RSlemma51}
F(\phi_1(z,w))=\begin{pmatrix}\bar A&\bar B\\\bar C&\bar D\end{pmatrix},
\quad
F(\phi_2(z,w))=\begin{pmatrix}D&C\\B&A\end{pmatrix},
\quad
F(\phi_3(z,w))=\begin{pmatrix}\bar D&\bar C\\\bar B&\bar A\end{pmatrix}
\end{equation}
and
\begin{equation}\label{eq:RSlemma52}
F(\phi_1(z,w))=\begin{pmatrix}\bar A&-\bar B\\-\bar C&\bar D\end{pmatrix}.
\end{equation}
\end{lemma}

Let $c_1:[0,1]\to M$ be a curve starting at $c_1(0)=(0,1)$ 
whose projection to the $z$-plane is an embedded curve in the first quadrant, 
and whose endpoint $c_1(1)$ has a $z$ coordinate so that $z\in\mathbb{R}$ 
and $1<z<a$.  
Let $c_2(t):[0,1]\to M$ be a curve starting at $c_2(0)=(0,1)$ 
whose projection to the $z$-plane is an embedded curve in the first quadrant, 
and whose endpoint $c_2(1)$ has a $z$ coordinate so that $z\in\mathbb{R}$ 
and $z>a$.  
With $F(0,1)=e_0$, we solve Equation \eqref{eq:F^-1dF} along these two paths 
to find 
\[ 
F(c_1(1))=\begin{pmatrix}A_1&B_1\\C_1&D_1\end{pmatrix},
\quad\text{and}\quad
F(c_2(1))=\begin{pmatrix}A_2&B_2\\C_2&D_2\end{pmatrix}.
\]
Let $\tau_j$ be the deck transformation of $\widetilde{M}$ associated to the 
homotopy class of $\gamma_j$ ($j=1,2,3$). 
\begin{itemize}
\item
Traveling about the loop $\gamma_1$, it follows from Equations 
\eqref{eq:RSlemma51} and \eqref{eq:RSlemma52} 
that $F\circ\tau_1=F\Phi_1$, where 
\[ 
\Phi_1 :=\begin{pmatrix}\bar{A}_1&-\bar{C}_1\\
                       -\bar{B}_1&\bar{D}_1\end{pmatrix}
         \begin{pmatrix}D_1&-C_1\\-B_1&A_1\end{pmatrix}
         \begin{pmatrix}\bar{D}_1&\bar{B}_1\\\bar{C}_1&\bar{A}_1\end{pmatrix}
         \begin{pmatrix}A_1&B_1\\C_1&D_1\end{pmatrix}.
\]
\item
Traveling about the loop $\gamma_2$, it follows from 
Equation \eqref{eq:RSlemma51} that $F\circ\tau_2=F\Phi_2$, where 
\[ 
\Phi_2 :=\begin{pmatrix}\bar{D}_2&-\bar{B}_2\\-\bar{C}_2&\bar{A}_2\end{pmatrix}
         \begin{pmatrix}A_2&B_2\\C_2&D_2\end{pmatrix}.
\]
\item
Traveling about $\gamma_3$, it follows from 
Equation \eqref{eq:RSlemma51} that $F\circ\tau_3=F\Phi_3$, where
\[ 
\Phi_3 :=\begin{pmatrix}\bar{A}_2&-\bar{C}_2\\
                       -\bar{B}_2&\bar{D}_2\end{pmatrix}
         \begin{pmatrix}D_2&C_2\\B_2&A_2\end{pmatrix}.
\]
\end{itemize}

We now wish to change the initial condition from $F(0,1)=e_0$ to 
\[
F(0,1)=P=\begin{pmatrix}P_{11}&P_{12}\\P_{21}&P_{22}\end{pmatrix}\in 
         SL(2,\mathbb{C})
\]
so that the $SU(1,1)$ conditions on all three loops $\gamma_1$, $\gamma_2$ and 
$\gamma_3$ will be solved. That is, we now find a constant matrix $P$ so that 
\[
P^{-1}\Phi_1P\quad\text{and}\quad P^{-1}\Phi_2P\quad\text{and}\quad 
P^{-1}\Phi_3P
\]
are all in $SU(1,1)$.  

To do this, we prepare several lemmas.
First of all, we show the following two lemmas about the loops $\gamma_2$ and 
$\gamma_3$:

\begin{lemma}\label{lm:1}
$\Phi_2$ and $\Phi_3$ can be written as follows:
\[
\Phi_2 =\begin{pmatrix}\psi_{11}&i\psi_{12}\\
                      i\psi_{21}&\bar\psi_{11}\end{pmatrix}, \qquad
\Phi_3 =\begin{pmatrix}\bar\psi_{11}&i\psi_{21}\\
                      i\psi_{12}&\psi_{11}\end{pmatrix},
\]
where $\psi_{11}\in\mathbb{C}$ and $\psi_{12},\psi_{21}\in\mathbb{R}$. 
\end{lemma}

\begin{proof}
By direct calculation and setting
\[
 \psi_{11}:=A_2\bar D_2-\bar B_2C_2, \quad
i\psi_{12}:=B_2\bar D_2-\bar B_2D_2, \quad
i\psi_{21}:=\bar A_2C_2-A_2\bar C_2,
\]
we get the conclusion. 
\end{proof}

Since $P\in SL(2,\mathbb{C})$, direct computation gives: 

\begin{lemma}\label{lm:2}
\begin{enumerate}
\item 
     For $P^{-1}\Phi_2P$ to be in $SU(1,1)$, we need 
     \begin{equation}\label{eq:psi}
     \left\{\begin{array}{l}
            (P_{12}P_{21}-\overline{P_{12}P_{21}})(\psi_{11}-\bar\psi_{11}) \\
           \hspace{74pt}
           -(P_{11}P_{12}-\overline{P_{11}P_{12}})i\psi_{21}
           +(P_{21}P_{22}-\overline{P_{21}P_{22}})i\psi_{12}=0, \\
            (P_{11}P_{21}-\overline{P_{12}P_{22}})(\psi_{11}-\bar\psi_{11})
           -(P_{11}^2-\bar P_{12}^2)i\psi_{21}
           +(P_{21}^2-\bar P_{22}^2)i\psi_{12}=0.
     \end{array}\right.
     \end{equation}
\item 
     For $P^{-1}\Phi_3P$ to be in $SU(1,1)$, we need 
     \begin{equation}\label{eq:tpsi}
     \left\{\begin{array}{l}
            (P_{12}P_{21}-\overline{P_{12}P_{21}})(\bar\psi_{11}-\psi_{11}) \\
           \hspace{74pt}
           +(P_{21}P_{22}-\overline{P_{21}P_{22}})i\psi_{21}
           -(P_{11}P_{12}-\overline{P_{11}P_{12}})i\psi_{12}=0, \\
            (P_{11}P_{21}-\overline{P_{12}P_{22}})(\bar\psi_{11}-\psi_{11})
           +(P_{21}^2-\bar P_{22}^2)i\psi_{21}
           -(P_{11}^2-\bar P_{12}^2)i\psi_{12}=0.
     \end{array}\right.
     \end{equation}
\end{enumerate}
\end{lemma}

If 
\begin{align}
P_{11}P_{12}-\overline{P_{11}P_{12}}&=P_{21}P_{22}-\overline{P_{21}P_{22}},
 \label{eq:ac-bdinR} \\
P_{11}^2-\bar P_{12}^2&=P_{21}^2-\bar P_{22}^2 \label{eq:c^2-a^2}
\end{align}
hold, then Equations \eqref{eq:psi} and \eqref{eq:tpsi} are equivalent. 
But we do not want both $P_{11}^2-\bar P_{12}^2$ and $P_{21}^2-\bar P_{22}^2$ 
to be zero unless 
$P_{11}P_{21}-\overline{P_{12}P_{22}}=0$. 

Next, we show the following two lemmas about the loop $\gamma_1$: 

\begin{lemma}\label{lm:3}
$\Phi_1$ can be written as follows:
\[
\Phi_1=\begin{pmatrix}\varphi_{11}&\varphi_{12}\\
                     -\bar\varphi_{12}&\varphi_{22}\end{pmatrix},
\]
where $\varphi_{11},\varphi_{22}\in\mathbb{R}$ and 
$\varphi_{12}\in\mathbb{C}$. 
\end{lemma}

\begin{proof}
By direct calculation and setting
\begin{align*}
\varphi_{11}&:=|\bar A_1D_1+B_1\bar C_1|^2-(\bar A_1C_1+A_1\bar C_1)^2, \\
\varphi_{22}&:=|\bar A_1D_1+B_1\bar C_1|^2-(\bar B_1D_1+B_1\bar D_1)^2, \\
\varphi_{12}&:=(\bar A_1D_1+B_1\bar C_1)
               (\bar B_1D_1+B_1\bar D_1-\bar A_1C_1-A_1\bar C_1),
\end{align*}
we get the conclusion. 
\end{proof}

Direct computation gives: 

\begin{lemma}\label{lm:4}
For $P^{-1}\Phi_1P$ to be in $SU(1,1)$, we need 
\begin{equation}\label{eq:phi}
\left\{\begin{array}{l}
       (\overline{P_{11}P_{22}}+P_{12}P_{21})\varphi_{11}
      -(P_{11}P_{22}+\overline{P_{12}P_{21}})\varphi_{22} \\
      \hspace{74pt}
      +(\overline{P_{11}P_{12}}+P_{21}P_{22})\varphi_{12}
      +(P_{11}P_{12}+\overline{P_{21}P_{22}})\bar\varphi_{12}=0, \\
       (P_{11}P_{21}+\overline{P_{12}P_{22}})(\varphi_{11}-\varphi_{22})
      +(\bar P_{12}^2+P_{21}^2)\varphi_{12}
      +(P_{11}^2+\bar P_{22}^2)\bar\varphi_{12}=0. 
\end{array}\right.
\end{equation}
\end{lemma}

\begin{remark}
Note that if we assume Equation \eqref{eq:c^2-a^2}, then the second equation 
of \eqref{eq:phi} can be replaced by 
\begin{equation}\label{eq:y+bary}
 (P_{11}P_{21}+\overline{P_{12}P_{22}})(\varphi_{11}-\varphi_{22})
+(P_{11}^2+\bar P_{22}^2)(\varphi_{12}+\bar\varphi_{12})=0. 
\end{equation}
\end{remark}

We set
\begin{equation}
              P=
P(\alpha,\beta)=\begin{pmatrix}P_{11}&P_{12}\\P_{21}&P_{22}\end{pmatrix}
               =\begin{pmatrix}\alpha&\varepsilon\beta\\
                               \alpha&-\varepsilon\beta\end{pmatrix},
\label{eq:initialP}
\end{equation}
where $\alpha,\beta\in\mathbb{C}$ satisfy $\alpha\beta =-\varepsilon /2$, 
and $\varepsilon$ is either $+1$ or $-1$. 
Then $\det P=1$ and Equations \eqref{eq:ac-bdinR} and \eqref{eq:c^2-a^2} hold, 
and hence Equations \eqref{eq:psi} and \eqref{eq:tpsi} are equivalent.  
Furthermore, we see that the first equations of both \eqref{eq:psi} and 
\eqref{eq:phi} vanish. 
Thus Equations \eqref{eq:psi} and \eqref{eq:tpsi} reduce 
\begin{equation}\label{eq:tobef2}
(\alpha^2+\bar\beta^2)(\psi_{11}-\bar\psi_{11})
+(\alpha^2-\bar\beta^2)i(\psi_{12}-\psi_{21})=0
\end{equation}
and Equations \eqref{eq:phi} reduce to 
\begin{equation}\label{eq:tobef1}
(\alpha^2-\bar\beta^2)(\varphi_{11}-\varphi_{22})
+(\alpha^2+\bar\beta^2)(\varphi_{12}+\bar\varphi_{12})=0. 
\end{equation}

\begin{theorem}
Let $(G,Q)=(w,cdzdw/w)$ be the Weierstrass data on $M$ defined as in 
\eqref{eq:riemsurf}. 
Let $F:\widetilde{M}\to SL(2,\mathbb{C})$ be the holomorphic null immersion 
so that $F$ satisfies \eqref{eq:F^-1dF} with initial condition 
$F(0,1)=P(\alpha, \beta)$ as in \eqref{eq:initialP}. 
We set 
\[ 
F(c_1(1))=\begin{pmatrix}A_1'&B_1'\\C_1'&D_1'\end{pmatrix},
\quad\text{and}\quad
F(c_2(1))=\begin{pmatrix}A_2'&B_2'\\C_2'&D_2'\end{pmatrix}.
\]
Then the following two conditions are equivalent:
\begin{enumerate}
\item $F$ satisfies the $SU(1,1)$ condition, 
\item $\alpha$ and $\beta$ satisfy 
\begin{equation}\label{eq:f1=f2}
\begin{split}
f_1:&=-\frac{\bar A_1'C_1'+A_1'\bar C_1'+\bar B_1'D_1'+B_1'\bar D_1'}
            {\bar A_1'D_1'+A_1'\bar D_1'+\bar B_1'C_1'+B_1'\bar C_1'} \\
    &=-\frac{\bar A_2'C_2'-A_2'\bar C_2'+\bar B_2'D_2'-B_2'\bar D_2'}
            {\bar A_2'D_2'-A_2'\bar D_2'+\bar B_2'C_2'-B_2'\bar C_2'}=:f_2
\end{split}
\end{equation}
and the absolute value of this number is greater than $1$. 
\end{enumerate}
\end{theorem}

\begin{proof}
By \eqref{eq:tobef2}, we have
\[
\varepsilon\frac{\bar\alpha^2+\beta^2}{\bar\alpha^2-\beta^2}
         =-\frac{i(\psi_{12}-\psi_{21})}{\psi_{11}-\bar\psi_{11}}
         =-\frac{\bar A_2'C_2'-A_2'\bar C_2'+\bar B_2'D_2'-B_2'\bar D_2'}
                {\bar A_2'D_2'-A_2'\bar D_2'+\bar B_2'C_2'-B_2'\bar C_2'}.
\]
Also, by \eqref{eq:tobef1}, we have
\[
\varepsilon\frac{\bar\alpha^2+\beta^2}{\bar\alpha^2-\beta^2}
         =-\frac{\varphi_{11}-\varphi_{22}}{\varphi_{12}+\bar\varphi_{12}}
         =-\frac{\bar A_1'C_1'+A_1'\bar C_1'+\bar B_1'D_1'+B_1'\bar D_1'}
                {\bar A_1'D_1'+A_1'\bar D_1'+\bar B_1'C_1'+B_1'\bar C_1'}.
\]
Moreover, since $\alpha =-\varepsilon /2\beta$, 
\[
 \varepsilon\frac{\bar\alpha^2+\beta^2}{\bar\alpha^2-\beta^2}
=\varepsilon\frac{1+4|\beta|^4}{1-4|\beta|^4}
\]
whose absolute value is greater than $1$ for any $\beta\in\mathbb{C}$, 
proving the theorem. 
\end{proof}

Therefore, if Equation \eqref{eq:f1=f2} holds, we choose $\alpha$ and $\beta$ 
and $\varepsilon$ so that 
\[
f_1=\varepsilon\frac{1+4|\beta|^4}{1-4|\beta|^4}=f_2
\]
and then the $SU(1,1)$ condition is satisfied. 

\begin{lemma}
If some $\alpha$, $\beta$ satisfy \eqref{eq:f1=f2}, we may assume 
$\alpha ,\beta\in\mathbb{R}$ and that \eqref{eq:f1=f2} still holds. 
\end{lemma}

\begin{proof}
Since $\alpha\beta =-\varepsilon /2$, there exists $r>0$ and 
$\theta\in [0,2\pi )$ so that 
\[
\alpha=re^{i\theta}\quad\text{and}\quad
\beta=\frac{-\varepsilon}{2r}e^{-i\theta}.
\]
Also, if $P^{-1}\Phi_jP\in SU(1,1)$ for $j=1,2,3$, then 
$(PU)^{-1}\Phi_j(PU)\in SU(1,1)$ for any $U\in SU(1,1)$ and $j=1,2,3$. 
Thus, setting $U=\text{diag}(e^{-i\theta},e^{i\theta})$, we see that
\[
PU=\begin{pmatrix}re^{i\theta}&(-1/2r)e^{-i\theta}\\
                  re^{i\theta}&(1/2r)e^{-i\theta}\end{pmatrix}
   \begin{pmatrix}e^{-i\theta}&0\\0&e^{i\theta}\end{pmatrix}
  =\begin{pmatrix}r&-1/2r\\r&1/2r\end{pmatrix}
\]
and hence each entry of $PU$ is real. 
\end{proof}

\begin{example}\label{ex:g1cat}
Now, in order to show 
the existence of a one-parameter family of weakly-complete CMC $1$ faces 
of genus $1$ with two ends which satisfy equality of the 
Osserman-type inequality,  we find values $c\in\mathbb{R}\setminus\{0\}$ and 
$a>1$ so that $|f_1|=|f_2|>1$ and $f_1=f_2$. 
By numerical experiments using Mathematica, we found such values 
(see Figure \ref{fg:a=2}).  
Also, by Corollary \ref{co:a=2} in Appendix \ref{app:ends}, 
we see that the ends are elliptic ends for $c<0$ 
(resp. hyperbolic ends for $c>0$).  
\end{example}

\begin{figure}[htbp] %%%%%%%%%%%%%%%%%%%%%%%%%%%%%%%%%%%%%%%%%%%%%%%%%%%
\begin{center}
 \includegraphics[width=.66\linewidth]{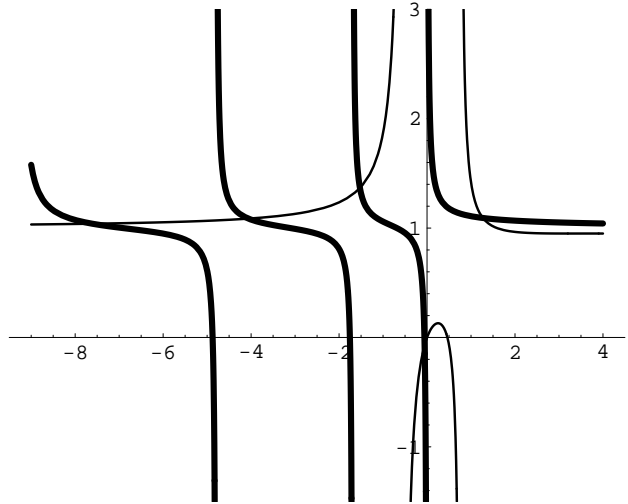}
\end{center}
\caption{The function $f_1$ (thin curve) and $f_2$ (thick curve) when $a=2$. 
         The horizontal axis represents $c$, and the vertical axis represents 
         $f_1$ and $f_2$.  We see that $f_1$ and $f_2$ intersect 6 times 
         for $c\in (-9,4)$, at $c\approx -7.6119$, $c\approx -4.06015$, 
         $c\approx -1.526035$, $c\approx -0.55$, $c\approx 1.26988$, 
         and $f_1=f_2>1$ except for $c\approx -0.55$.}
\label{fg:a=2}
\end{figure} %%%%%%%%%%%%%%%%%%%%%%%%%%%%%%%%%%%%%%%%%%%%%%%%%%%%%%%%%%%

\begin{figure}[htbp] %%%%%%%%%%%%%%%%%%%%%%%%%%%%%%%%%%%%%%%%%%%%%%%%%%%
\begin{center}
\begin{tabular}{cc}
 \includegraphics[width=.50\linewidth]{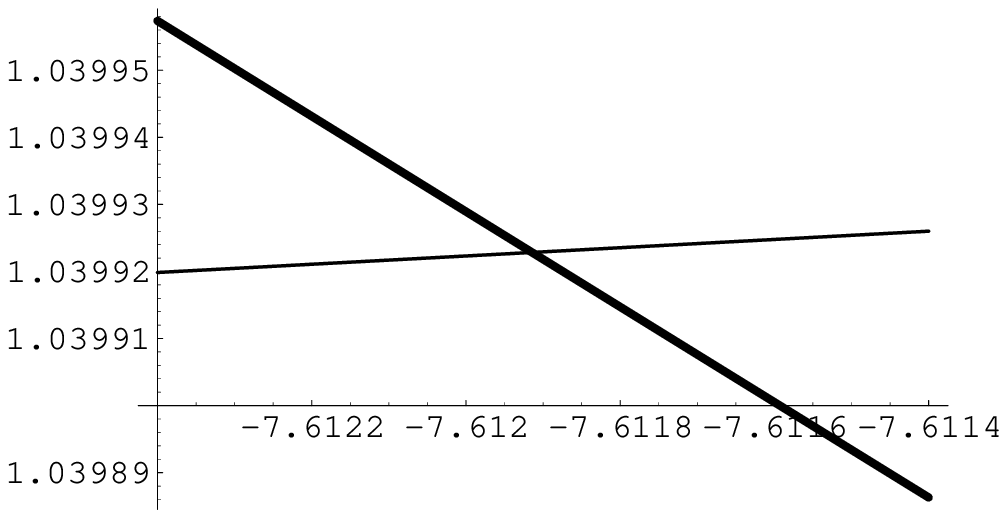} &
 \includegraphics[width=.40\linewidth]{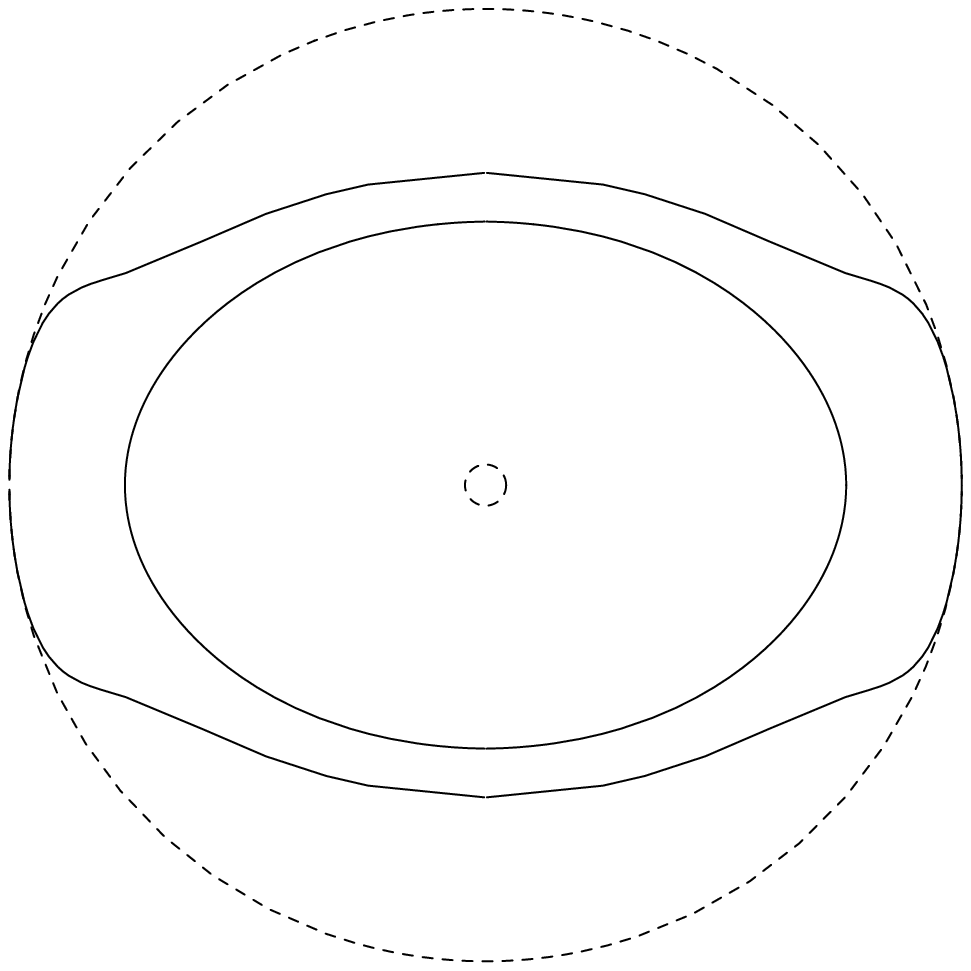} \\
 \includegraphics[width=.50\linewidth]{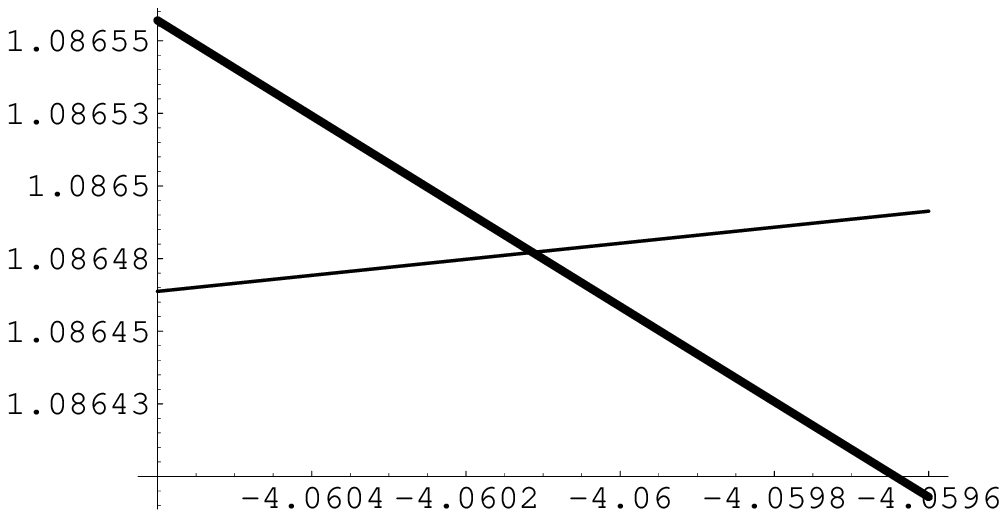} &
 \includegraphics[width=.40\linewidth]{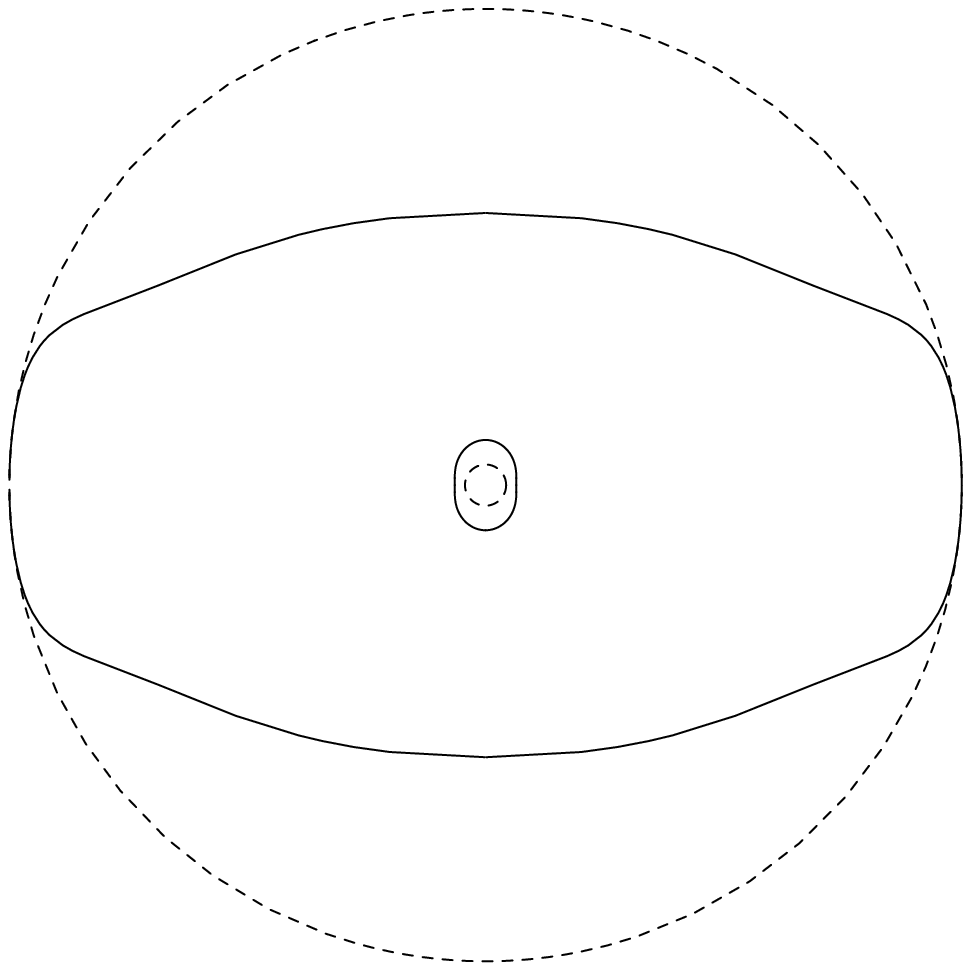} \\
 \includegraphics[width=.50\linewidth]{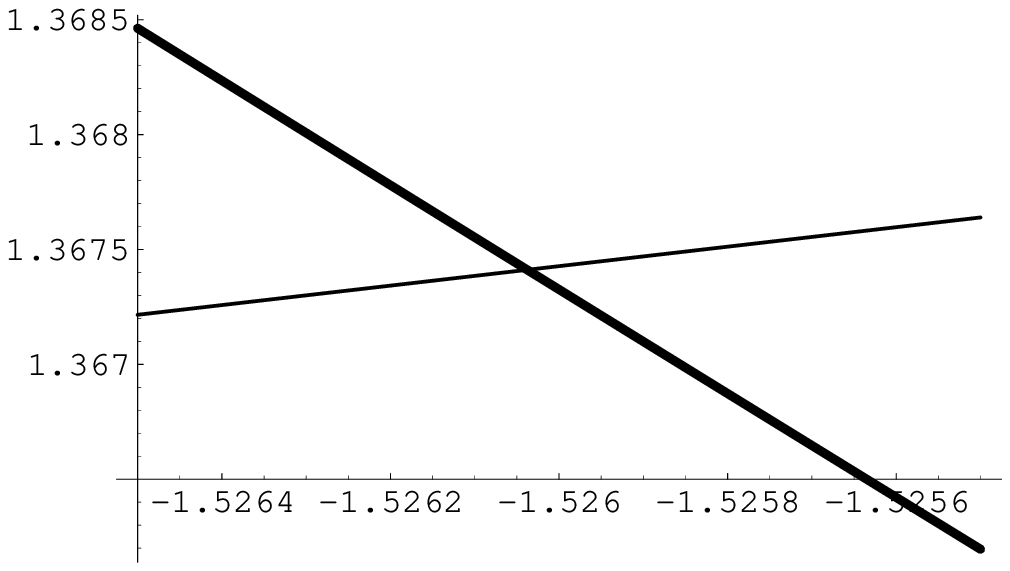} &
 \includegraphics[width=.40\linewidth]{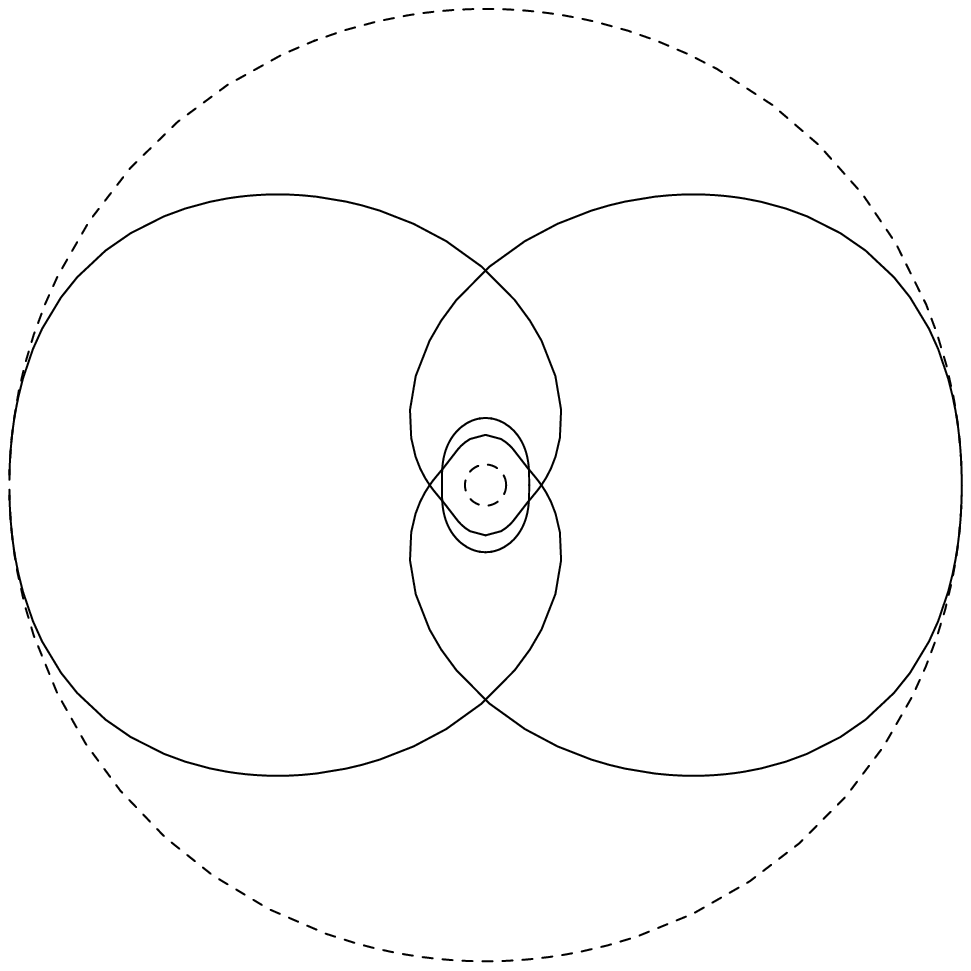} 
\end{tabular}
\end{center}
\caption{Left: The function $f_1$ (thin curve) and $f_2$ (thick curve) when 
         $a=2$.  The horizontal axis represents $c$, and the vertical axis 
         represents $f_1$ and $f_2$.  We see that $f_1,f_2>1$ for 
         $c\in (-7.6124,-7.6114)$ in the first row, 
         $c\in (-4.0606,-4.0596)$ in the second row and 
         $c\in (-1.5265,-1.5255)$ in the third row, and $f_1=f_2$ at some such 
         value of $c$ in each case, and $a=2>1$.  
         Right: Symmetry curves in the CMC $1$ face in Example \ref{ex:g1cat}
         intersect the plane $\{(y_1,y_2,y_3)\in\mathscr{H}\,|\,y_2=0\}$, 
         with $a=2$ and $c=-7.6119$ (resp. $c=-4.06015$, $c=-1.526035$).}
\label{fg:g1cat-0}
\end{figure} %%%%%%%%%%%%%%%%%%%%%%%%%%%%%%%%%%%%%%%%%%%%%%%%%%%%%%%%%%%

\begin{figure}[htbp] %%%%%%%%%%%%%%%%%%%%%%%%%%%%%%%%%%%%%%%%%%%%%%%%%%%
\begin{center}
 \includegraphics[width=.40\linewidth]{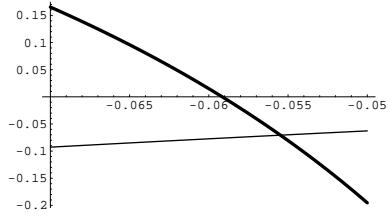}
\end{center}
\caption{The function $f_1$ (thin curve) and $f_2$ (thick curve) when $a=2$. 
         The horizontal axis represents $c$, and the vertical axis represents 
         $f_1$ and $f_2$.  We see that $f_1=f_2$ at some value of 
         $c\in (-0.07,0.05)$ but $|f_1|=|f_2|<1$ at this value of $c$.}
\label{fg:0}
\end{figure} %%%%%%%%%%%%%%%%%%%%%%%%%%%%%%%%%%%%%%%%%%%%%%%%%%%%%%%%%%%

\begin{figure}[htbp] %%%%%%%%%%%%%%%%%%%%%%%%%%%%%%%%%%%%%%%%%%%%%%%%%%%
\begin{center}
\begin{tabular}{cc}
 \includegraphics[width=.50\linewidth]{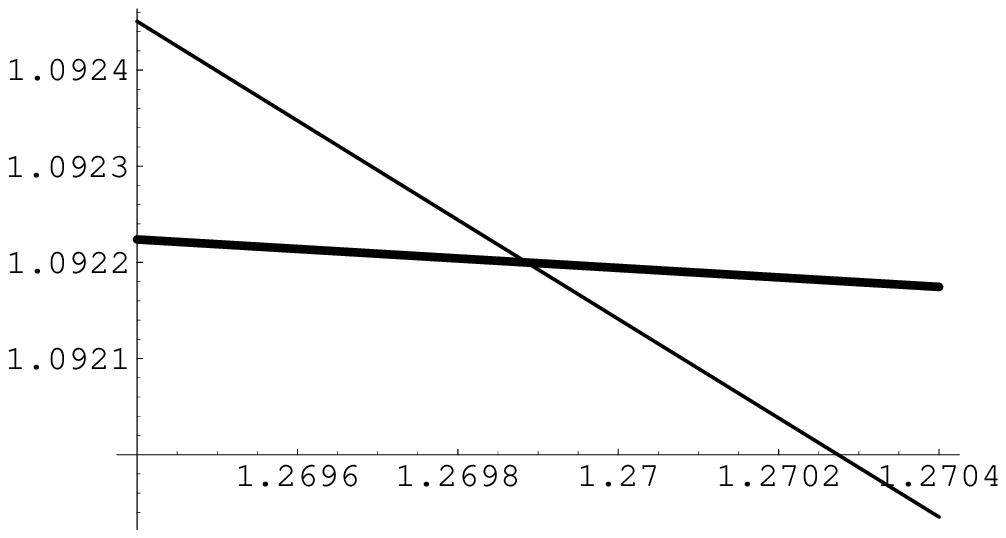} &
 \includegraphics[width=.40\linewidth]{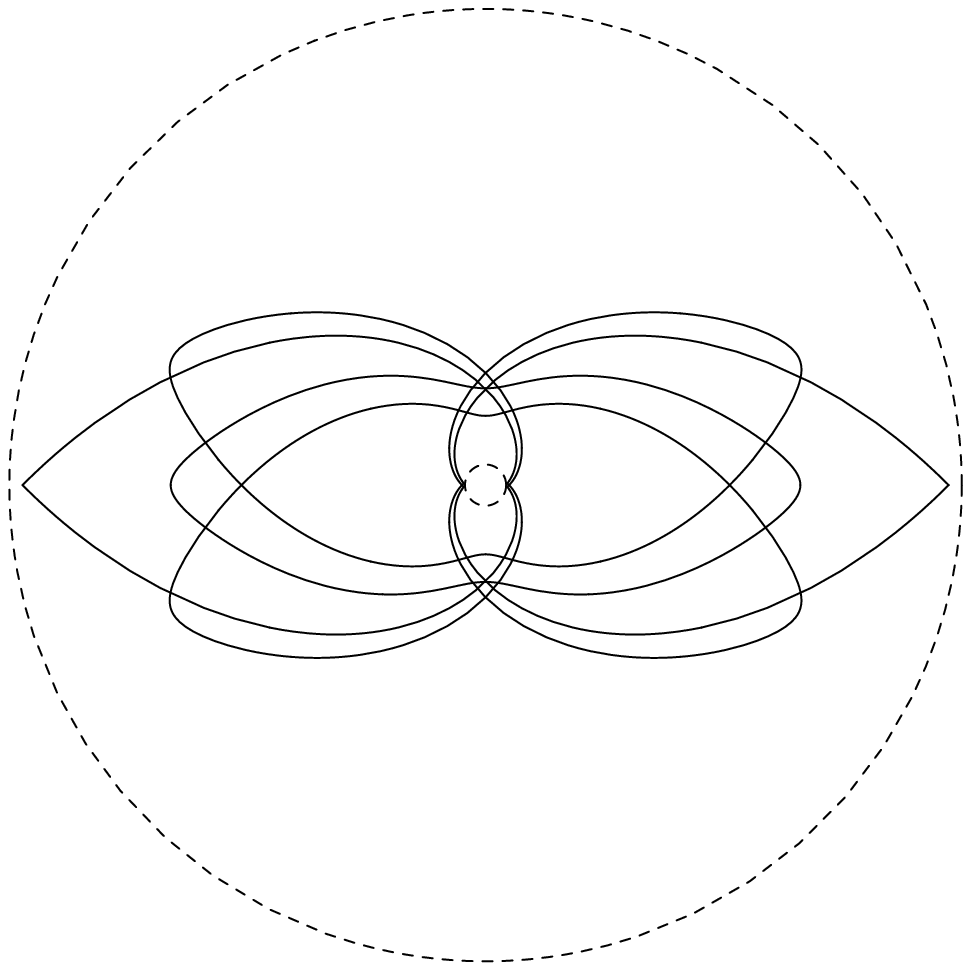} 
\end{tabular}
\end{center}
\caption{Left: The function $f_1$ (thin curve) and $f_2$ (thick curve) when 
         $a=2$.  The horizontal axis represents $c$, and the vertical axis 
         represents $f_1$ and $f_2$.  We see that $f_1,f_2>1$ for 
         $c\in (1.2694,1.2704)$, and $f_1=f_2$ at some such value of $c$, 
         and $a=2>1$.
         Right: Symmetry curves in the CMC $1$ face in Example \ref{ex:g1cat}
         intersect the plane $\{(y_1,y_2,y_3)\in\mathscr{H}\,|\,y_2=0\}$, 
         with $a=2$ and $c=1.26988$.}
\label{fg:g1cat+1}
\end{figure} %%%%%%%%%%%%%%%%%%%%%%%%%%%%%%%%%%%%%%%%%%%%%%%%%%%%%%%%%%%

\appendix %%%%%%%%%%%%%%%%%%%%%%%%%%%%%%%%%%%%%%%%%%%%%%%%%%%%%%%%%%%%%%

\section{Criteria for the types of ends} %%%%%%%%%%%%%%%%%%%%%%%%%%%%%%%
\label{app:ends} %%%%%%%%%%%%%%%%%%%%%%%%%%%%%%%%%%%%%%%%%%%%%%%%%%%%%%%

Here we give the criteria for when an end of the genus $1$ catenoid 
given in \eqref{eq:riemsurf} and \eqref{eq:W-data} is elliptic or 
hyperbolic. 
First we give the following lemma:

\begin{lemma}\label{lm:indep}
Let $f:M\to\mathbb{S}^3_1$ be a CMC $1$ face 
and $\gamma$ a loop in $M$.  Then the eigenvalues of the monodromy matrix 
with respect to $\gamma$ do not depend on the choice of the holomorphic null 
lift $F$ of $f$. 
\end{lemma}

\begin{proof}
Let $(G,Q)$ be a Weierstrass data of $f$ and 
$F_1, F_2:\widetilde{M}\to SL(2,\mathbb{C})$ solutions of Equation 
\eqref{eq:F^-1dF}. 
Then there exists a constant $B\in SL(2,\mathbb{C})$ such that 
$F_1=F_2B$.
Let $\tau$ be the deck transformation of $\widetilde{M}$ associated to the 
homotopy class of $\gamma$ and $\Phi_j$ ($j=1,2$) the monodromy matrix of 
$F_j$ with respect to $\gamma$.  Then 
\begin{align*}
F_1\circ\tau &= F_1\Phi_1=F_2B\Phi_1 \\
             &= (F_2B)\circ\tau =F_2\Phi_2B.
\end{align*}
Thus $\Phi_1=B^{-1}\Phi_2B$ and hence the eigenvalues of $\Phi_1$ and 
$\Phi_2$ are the same, proving the lemma.  
\end{proof}

So to determine the type of an end, 
we can take any holomorphic null lift $F$ of $f$. 
Let $F=(F_{jk})_{j,k=1,2}:\widetilde{M}\to SL(2,\mathbb{C})$ be a holomorphic 
null lift of $f$.  Direct calculation shows that 
\begin{align}
\frac{d^2F_{1j}}{dz^2}-\frac{1}{w}\frac{dw}{dz}\frac{dF_{1j}}{dz}
                      +\frac{c}{w}\frac{dw}{dz}F_{1j} &=0, 
\tag{E.1}\label{eq:E1} \\
\frac{d^2F_{2j}}{dz^2}+\frac{1}{w}\frac{dw}{dz}\frac{dF_{2j}}{dz}
                      +\frac{c}{w}\frac{dw}{dz}F_{2j} &=0 
\tag{E.2}\label{eq:E2}
\end{align} 
for $j=1,2$.  
We consider the end $(z,w)=(\infty,1)$.  
Let $\Delta^*\subset M$ be a neighborhood of $(z,w)=(\infty,1)$. 
We set $\zeta =1/z$.  Without loss of generality we may assume 
$\Delta^*=\{\zeta\in\mathbb{C}\,|\,0<|\zeta|<1\}$. 
We set $\Delta =\Delta^*\cup\{0\}$.  
Then Equations \eqref{eq:E1} and \eqref{eq:E2} become 
\begin{align}
\zeta^2\frac{d^2F_{1j}}{d\zeta^2}+\zeta p_1(\zeta)\frac{dF_{1j}}{dz}
                                 +q(\zeta) F_{1j} &=0, 
\label{eq:E1zeta} \\
\zeta^2\frac{d^2F_{2j}}{dz^2}+\zeta p_2(\zeta)\frac{dF_{2j}}{dz}
                             +q(\zeta)F_{2j} &=0, 
\label{eq:E2zeta}
\end{align} 
where
\[
p_1(\zeta)=2-\frac{1}{w}\frac{dw}{d\zeta}\zeta, \qquad
p_2(\zeta)=2+\frac{1}{w}\frac{dw}{d\zeta}\zeta \qquad\text{and}\qquad
q  (\zeta)= -\frac{c}{w}\frac{dw}{d\zeta}.
\]
Note that 
\[
\frac{1}{w}\frac{dw}{d\zeta}
=\frac{(1-a)(a\zeta^2+1)}{(\zeta^2-1)(a^2\zeta^2-1)}
=(1-a)+\mathcal{O}(\zeta^2).
\]
Fundamental systems of solutions $\{X_1,X_2\}$ of 
Equation \eqref{eq:E1zeta} and $\{Y_1,Y_2\}$ of Equation \eqref{eq:E2zeta} 
can be chosen as 
\begin{align}
X_1=\zeta^{(-1+m)/2}\xi_1(\zeta),\qquad &
X_2=\zeta^{(-1-m)/2}\xi_2(\zeta)+k_1 X_1\log\zeta, 
\label{eq:X_j}\\
Y_1=\zeta^{(-1+m)/2}\eta_1(\zeta),\qquad &
Y_2=\zeta^{(-1-m)/2}\eta_2(\zeta)+k_2 Y_1\log\zeta,
\label{eq:Y_j}
\end{align}
where
\[
m=\sqrt{1-4c(a-1)},
\]
and $\xi_j$ and $\eta_j$ ($j=1,2$) are holomorphic functions on $\Delta$ with 
$\xi_j(0)\ne 0$ and $\eta_j(0)\ne 0$, and the constant $k_1$ (resp. $k_2$) is 
called the log-term coefficient of the solutions of Equation \eqref{eq:E1zeta} 
(resp. \eqref{eq:E2zeta}).  See, for example, \cite[Appendix A]{RUY4}. 

Although $m$ can be either a positive real or is purely imaginary, 
here we only consider the case $m\not\in\mathbb{Z}$.  
In this case, it is known that 
\[
k_1=k_2=0.
\]
Moreover, we have the following lemma:

\begin{lemma}\label{lm:UY1-53}
There exists a matrix $\Lambda\in SL(2,\mathbb{C})$ such that 
\begin{equation}
\label{eq:FLambda}
F\Lambda =\begin{pmatrix}
          \zeta^{(-1+m)/2}A(\zeta) & \zeta^{(-1-m)/2}B(\zeta) \\
          \zeta^{(-1+m)/2}C(\zeta) & \zeta^{(-1-m)/2}D(\zeta) 
          \end{pmatrix},
\end{equation}
where $A$, $B$, $C$ and $D$ are holomorphic functions on $\Delta$ such that 
$A(0)$, $B(0)$, $C(0)$ and $D(0)$ are all nonzero.  
\end{lemma}

\begin{proof}
Since $f$ is not totally umbilic, $F_{11}$ and $F_{12}$ are linearly 
independent and are linear combinations of the $X_1$ and $X_2$ in Equation 
\eqref{eq:X_j}.  Then there exists a matrix 
$\Lambda\in SL(2,\mathbb{C})$ such that 
\[
F\Lambda 
 = \begin{pmatrix}
   \zeta^{(-1+m)/2}A(\zeta) & \zeta^{(-1-m)/2}B(\zeta) \\
   C_1\zeta^{(-1+m)/2}\eta_1(\zeta)+C_2\zeta^{(-1-m)/2}\eta_2(\zeta) & 
   D_1\zeta^{(-1+m)/2}\eta_1(\zeta)+D_2\zeta^{(-1-m)/2}\eta_2(\zeta)
   \end{pmatrix},
\]
where $C_j$ and $D_j$ ($j=1,2$) are constants, and $A$ and $B$ are holomorphic 
functions on $\Delta$ such that $A(0)\ne 0$ and $B(0)\ne 0$. 
Since $F\Lambda\in SL(2,\mathbb{C})$, we have 
\begin{align*}
1 &=\det (F\Lambda) \\
  &=D_1A(\zeta)\eta_1(\zeta)\zeta^{-1+m}
    +\bigl(D_2A(\zeta)\eta_2(\zeta)-C_1B(\zeta)\eta_1(\zeta)\bigr)\zeta^{-1}
    -C_2B(\zeta)\eta_2(\zeta)\zeta^{-1-m}. 
\end{align*}
Since $A(0)$, $B(0)$ and $\eta_j(0)$ are all nonzero, 
it follows that $D_1=C_2=0$.  Setting $C(\zeta)=C_1\eta_1(\zeta)$ and 
$D(\zeta)=D_2\eta_2(\zeta)$, and noting that $C(0)\ne 0$ and $D(0)\ne 0$, 
we have the conclusion. 
\end{proof}

\begin{proposition}\label{pr:ends}
Let $M$ be the Riemann surface defined as in Equation \eqref{eq:riemsurf} and 
$f:M\to\mathbb{S}^3_1$ the CMC $1$ face constructed from the Weierstrass data 
as in Equation \eqref{eq:W-data}.  
Then the monodromy of an end is elliptic {\rm (}resp. hyperbolic{\rm )} 
if $m\in\mathbb{R}^+\setminus\mathbb{N}$ 
{\rm (}resp. $m\in i\mathbb{R}\setminus\{0\}${\rm )}. 
\end{proposition}

\begin{proof}
By Lemma \ref{lm:indep}, we can choose $F\Lambda$ as in Equation 
\eqref{eq:FLambda} as a holomorphic null lift of $f$. 
Let $\gamma$ be a loop around an end and $\tau$ the deck transformation of 
$\widetilde{M}$ associated to the homotopy class of $\gamma$.  
Then 
\[
(F\Lambda)\circ\tau
=(F\Lambda)\begin{pmatrix}-e^{m\pi i}&0\\0&-e^{m\pi i}\end{pmatrix}
=(F\Lambda)\Phi_\gamma . 
\]
Thus the eigenvalues of $\Phi_\gamma$ are $-e^{\pm m\pi i}$, 
which are in $\mathbb{S}^1$ (resp. $\mathbb{R}\setminus\{1\}$) if 
$m\in\mathbb{R}^+\setminus\mathbb{N}$ 
(resp. $m\in i\mathbb{R}\setminus\{0\}$), proving the proposition. 
\end{proof}

\begin{corollary}\label{co:a=2}
If $a=2$, then the monodromy of an end is elliptic 
{\rm (}resp. hyperbolic{\rm )} 
if $c<0$ {\rm (}resp. $c>0${\rm )}. 
\end{corollary}

 %%%%%%%%%%%%%%%%%%%%%%%%%%%%%%%%%%%%%%%%%%%%%%%%%%
\end{document}